\documentclass[12pt,reqno]{amsart}

\usepackage{amsthm,amsmath,amssymb,amsfonts,amscd,verbatim,latexsym}

\newtheorem{Theorem}{Theorem}[section]
\newtheorem{Lemma}[Theorem]{Lemma}
\newtheorem{Proposition}[Theorem]{Proposition}

\newtheorem{Example}[Theorem]{Example}
\newtheorem{Conjecture}[Theorem]{Conjecture}

\newcommand{\QQ}{\mathbf Q}
\newcommand{\ux}{\mathbf x} 
\newcommand{\A}{\mathcal A}
\newcommand{\C}{\mathcal C}
\newcommand{\Y}{\mathcal Y}
\newcommand{\Z}{\mathcal Z}
\newcommand{\ga}{\mathfrak a}
\newcommand{\gA}{\mathfrak A}
\newcommand{\gf}{\mathfrak f}

\newcommand{\tgf}{\widetilde{\mathfrak f}}
\renewcommand{\SS}{\mathfrak S}
\newcommand{\Sym}{\text{Sym}}
\newcommand{\spec}{\text{spec}\,}
\newcommand{\im}{\text{im}\,}
\renewcommand{\P}{\mathbf P}
\newcommand{\U}{\mathcal U}
\renewcommand{\O}{\mathcal O}
\renewcommand{\H}{\mathbb H} 
\newcommand{\HH}{\mathcal H} 
 
\newcommand{\K}{\mathcal K}
\newcommand{\F}{\mathbb F} 
 
\newcommand{\I}{\mathcal I} 
 
\newcommand{\Gordan}{\mathfrak G}

\renewcommand{\le}{\leqslant}
\renewcommand{\ge}{\geqslant}
\newcommand{\ra}{\rightarrow}

\newcommand{\lra}{\longrightarrow}

\renewcommand{\ge}{\geqslant} 

\newcommand{\demo}{\noindent {\sc Proof.}\;}

\begin{document}
\title[saturation sequence]
{On the saturation sequence of the rational normal curve} 
\author{Jaydeep Chipalkatti} 
\maketitle 

\bigskip 

\parbox{11.8cm}{ \small
{\sc Abstract:} Let $C \subseteq \P^d$ denote the rational normal curve of order $d$. Its 
homogeneous defining ideal $I_C \subseteq \QQ[a_0,\dots,a_d]$ admits an $SL_2$-stable filtration 
$J_2 \subseteq J_4 \subseteq \dots \subseteq I_C$ by sub-ideals such that the saturation of 
each $J_{2q}$ equals $I_C$. Hence, one can associate to $d$ a sequence of integers 
$(\alpha_1,\alpha_2,\dots)$ which encodes the degrees in which the successive inclusions 
in this filtration become trivial. In this paper we establish several lower and upper 
bounds on the $\alpha_q$, using \emph{inter alia} the methods of classical invariant theory.

\medskip 


\medskip 

} 

\medskip 

\thispagestyle{empty} 

\section{Introduction}
\subsection{} 
The rational normal curve of order $d$ in $\P^d$ and its homogeneous 
defining ideal usually make an obligatory appearance in textbooks on algebraic 
geometry\footnote{For instances, see~\cite[Exer.~A2.10]{Eisenbud}, \cite[Lecture 1]{Harris}, 
\cite[Ch.~IV, Exer.~{3.4}]{Hartshorne}.}. 
This is not without its reasons. The latter admits a winsome description as the ideal of 
maximal minors of a $2 \times d$ matrix of variables, usually called the catalecticant matrix (see~\S\ref{section.cat} below). 

However, this formulation disguises the fact that the ideal carries a nontrivial filtration 
which is invariant under the automorphisms of $\P^d$ fixing the curve. 
The object of this paper is to initiate a study of 
this filtration; the main results are described in~\S\ref{section.results} after the required 
notation is available. 

Throughout, the base field will be $\QQ$ (the field of rational numbers). 
Classical treatments of the necessary background in invariant theory may be found in 
\cite{GrYo,Salmon}, and more modern treatments 
in~\cite{Dolgachev,KungRota,Olver,Procesi,Sturmfels}. 
\subsection{Transvectants} 
Let $A(x_1,x_2)$ and $B(x_2,x_2)$ denote binary forms of orders $p,q$ respectively in the variables 
$\ux = \{x_1,x_2\}$. Their $r$-th 
transvectant\footnote{Usually $r$ is called the index of transvection.} is defined by the formula 
\begin{equation} 
(A,B)_r = \frac{(p-r)! \, (q-r)!}{p! \, q!} \, \sum\limits_{i=0}^r  
\, (-1)^i \, \binom{r}{i} \,
\frac{\partial^r A}{\partial x_1^{r-i} \, \partial x_2^i} \,
\frac{\partial^r B}{\partial x_1^i \, \partial x_2^{r-i}};   
\label{formula.transvectants} \end{equation} 
for $0 \le r \le \min(p,q)$. It is of order $p+q-2r$ in $\ux$. If $r > \min(p,q)$, then 
$(A,B)_r=0$. Moreover, $(A,B)_r = (-1)^r (B,A)_r$, and hence $(A,A)_r$ vanishes 
for odd values of $r$. 
\subsection{Representations of $SL_2$} 
For $p \ge 0$, let $S_p$ denote the set of binary forms of order $p$ in $\ux$ 
(with coefficients in $\QQ$). The 
group $SL_2 \, \QQ$ acts on $S_p$ as follows: for $g = \left( \begin{array}{cc} 
\gamma_{11} & \gamma_{12} \\ \gamma_{21} & \gamma_{22} \end{array} \right) 
\in SL_2$, 
\[ A(x_1,x_2) \stackrel{g}{\lra} 
A(\gamma_{11} \, x_1 + \gamma_{12} \, x_2, \gamma_{21} \, x_1 + 
\gamma_{22} \, x_2). \] 
Up to isomorphism, $\{S_p: p \ge 0\}$ is the set of all the finite-dimensional 
irreducible representations of $SL_2$, and each such representation splits as a direct sum of 
irreducibles (see~\cite[Ch.~10]{Procesi}). For any $p,q \ge 0$, there 
is a decomposition
\[ S_p \otimes S_q \simeq \bigoplus\limits_{r=0}^{\min(p,q)} \, S_{p+q-2r}, \] 
and the image of $A \otimes B$ via the projection map 
$S_p \otimes S_q \lra S_{p+q-2r}$, is the transvectant $(A,B)_r$. There is an isomorphism 
of $S_p$ with its dual representation $S_p^* = \text{Hom}(S_p,\QQ)$, which associates 
$A \in S_p$ with the functional $B \lra (A,B)_p$. 
\subsection{The ring of covariants} 
Fix an integer $d \ge 1$, and introduce variables $a_0,\dots,a_d$. 
Define the bigraded polynomial ring 
\[ \C = \QQ[a_0,\dots,a_d; x_1,x_2] = \bigoplus\limits_{m, n \ge 0} \, \C_{m,n}, \] 
where $m$ (respectively $n$) denotes the degree in the $a$-variables 
(respectively $\ux$-variables). Let 
\begin{equation} \F = \sum\limits_{i=0}^d \, \binom{d}{i} \, a_i \, x_1^{d-i} \, x_2^i \in \C, 
\label{F.gen} \end{equation}
denote the generic binary $d$-ic, and define $\A$ to be the 
smallest $\QQ$-subalgebra of $\C$ satisfying the following two properties: 
\begin{itemize} 
\item $\F \in \A$, 
\item if $T,T' \in \A$ are bihomogeneous elements, then $(T,T')_r \in \A$ for all $r \ge 0$. 
\end{itemize} 
In other words, $\A$ is spanned as a $\QQ$-vector space by all compound transvectant expressions 
\[ (\F,\F)_2, \; (\F,(\F,\F)_2)_5, \;  ((\F,\F)_2,(\F,\F)_4)_3, \dots \text{etc.}\] 
We have a bigraded decomposition, 
\[ \A = \bigoplus_{m,n} \, \A_{m,n}, \quad \text{where $\A_{m,n} = \C_{m,n} \cap \A$}. \] 
In classical literature $\A$ is called the ring of covariants\footnote{It is more common to 
define it as the invariant subring $\C^{SL_2}$, but our definition is equivalent.} 
(of a binary $d$-ic); and an element $\Phi \in \A_{m,n}$ is called a 
covariant of degree $m$ and order $n$. E.g., 
$(\F,(\F,\F)_2)_5$ is a covariant of degree $3$ and order $3d-14$. A covariant of order 
zero is called an invariant. 

It is a fundamental result due to Gordan that $\A$ is finitely generated as a $\QQ$-algebra 
(see~\cite[Ch.~VI]{GrYo}). E.g., if $d=4$, then $\A$ is generated by the elements 
\[ \F, \; (\F,\F)_2, \;  (\F,\F)_4, \;  (\F,(\F,\F)_2)_1, \;  (\F,(\F,\F)_2)_4; \] 
of degree-orders $(1,4),(2,4),(2,0),(3,6),(3,0)$ respectively. 
\subsection{} Now identify the generic form $\F$ with the natural trace element 
in $S_d^* \otimes S_d \simeq S_d \otimes S_d$; 
this amounts to letting $a_i = \frac{1}{d!} \, (-x_1)^i x_2^{d-i} \in S_d$. Then 
$R = \QQ[a_0,\dots,a_d]$ is identified with the symmetric algebra 
$\bigoplus\limits_{m \ge 0} \, \Sym^m \, S_d$. 
Consider the decomposition 
\[ R_m \simeq \Sym^m \, S_d \simeq \bigoplus\limits_n \, \, (S_n \otimes \QQ^{\eta_{m,n}}). \] 
A covariant 
$\Phi = \varphi_0 \, x_1^n + \varphi_1 \, x_1^{d-1} \, x_2  + \dots + \varphi_n \, x_2^n$
of degree-order $(m,n)$ gives an $SL_2$-equivariant morphism 
\[ S_n \lra R_m, \quad A \lra (A,\Phi)_n; \] 
and conversely, every such morphism arises from a covariant. Hence 
\[ \dim \A_{m,n} = \eta_{m,n} = \dim \text{Hom}_{SL_2}(S_n,R_m). 
\] E.g., for $d=6$, there is a decomposition 
\[ R_3  \simeq \Sym^3 S_6 \simeq S_{18} \oplus S_{14} \oplus S_{12} 
\oplus S_{10} \oplus S_8 \oplus (S_6 \otimes \QQ^2) \oplus S_2; \] 
in particular, $\dim \A_{3,6} = 2$. It 
is easy to verify that 
\[ \{(\F,(\F,\F)_2)_4, (\F,(\F,\F)_4)_2\} \] 
is a basis of $\A_{3,6}$. By contrast, since $\A_{3,8}$ is one-dimensional, 
the forms $(\F,(\F,\F)_2)_3$ 
and $(\F,(\F,\F)_4)_1$ must be dependent; in fact there is an identical relation 
$7 \, (\F,(\F,\F)_2)_3 - (\F,(\F,\F)_4)_1=0$. Such calculations in $\A$ can be 
carried out by using the classical symbolic calculus (see~\cite{GrYo}). 
\subsection{Quadratic covariants} 
Now let $e_d = [\frac{d}{2}]$, and write 
\[ \H_{2q} = (\F,\F)_{2q}, \qquad \text{for $1 \le q \le e_d$,} \] 
which is a covariant of degree $2$ and order $2d-4q$. (Usually $\H_2$ is called the Hessian of $\F$.) 
We have a decomposition 
\[ R_2 \simeq \Sym^2 \, S_d \simeq \bigoplus\limits_{q=0}^{e_d} \, S_{2d-4q},  \] 
in which the summand $S_{2d-4q}$ corresponds to the span of the coefficients of $\H_{2q}$. 
Define $W_{2q}$ to be the subspace of $R_2$ generated by all 
the coefficients of $\H_2,\H_4, \dots, \H_{2q}$, and let $J_{2q}$ be the ideal in $R$ generated 
by $W_{2q}$. This defines a filtration 
\begin{equation} 
J_2 \subsetneq J_4 \subsetneq \dots \subsetneq J_{2e_d}, 
\end{equation} 
which is nontrivial for all $d \ge 4$. 
\subsection{} \label{section.cat}
Now $\P S_d = \text{Proj} \, R$ is the space of binary $d$-ics 
(distinguished up to scalars). It is a classical result 
(see \cite[Proposition 2.23]{Olver}) that the following conditions are equivalent for $A \in S_d$. 
\begin{enumerate} 
\item $(A,A)_2=0$. 
\item $(A,A)_2=(A,A)_4=\dots =(A,A)_{2e_d}=0$. 
\item There exists a linear form $t_1 \, x_1 + t_2 \, x_2$, such that 
$A = (t_1 \, x_1 + t_2 \, x_2)^d$. 
\end{enumerate} 
It follows that the variety cut out by the ideal $J_{2e_d}$ is the rational normal curve 
$C = \{ [(t_1 \, x_1 + t_2 \, x_2)^d \, ]: t_1,t_2 \in \QQ\} \subseteq \P S_d$. Since 
the defining ideal $I_C \subseteq R$ is $SL_2$-stable and generated by quadrics, in fact 
$J_{2e_d} = I_C$. It may also be described as the ideal of maximal minors of the 
catalecticant matrix 
\[ \left[ \begin{array}{ccccc} 
a_0 & a_1 & \dots & a_{d-2} & a_{d-1} \\ 
a_1 & a_2 & \dots & a_{d-1} & a_d \end{array} \right]. \] 

The equivalence (1)$\iff$(3) implies that $J_2$ defines $C$ set-theoretically, but in fact a 
stronger statement holds. 
\begin{Proposition} \sl 
The saturation of $J_2$ equals $I_C$. 
\label{prop.saturation}\end{Proposition} 
\demo 
See \cite[Lemma~3.1]{AF}, as well as \S\ref{section.uv} below. \qed 

\subsection{} 
It follows that all the ideals $J_{2q}$ coincide in sufficiently high degrees. For 
$1 \le q \le e_d-1$, define 
\[ \alpha_q = \min \, \{m: (J_{2q})_t = (J_{2q+2})_t \; \; \text{for all $t \ge m$}\},  \] 
then $(\alpha_1,\dots,\alpha_{e_d-1})$ 
will be called the {\bf saturation sequence} of $d$. 
I am enclosing the table of saturation sequences 
for $d \le 20$. It was calculated in {\sc Macaulay-2}. 

\[ \begin{array}{rl} 
d & \text{saturation sequence} \\ \hline 
4 & (3) \\ 
5 & (3) \\ 
6 & (5,3) \\ 
7 & (4,3) \\ 
8 & (5,3,3) \\ 
9 & (5,3,3) \\ 
10 & (5,3,3,3) \\ 
11 & (5,3,3,3) \\ 
12 & (7,5,3,3,3) \\ 
13 & (5,4,3,3,3) \\ 
14 & (7,5,3,3,3,3) \\ 
15 & (6,5,3,3,3,3) \\ 
16 & (7,5,4,3,3,3,3) \\ 
17 & (7,5,4,3,3,3,3) \\ 
18 & (7,5,5,3,3,3,3,3) \\ 
19 & (7,5,4,3,3,3,3,3) \\ 
20 & (8,5,5,4,3,3,3,3,3) 
\end{array} \] 
Recall that the 
{\sl satiety} of $J_{2q}$ is defined to be the integer (cf.~\cite[p.~593]{BG})
\[ \min\{ \, m: (J_{2q})_t = (I_C)_t \; \; \text{for all $t \ge m$} \}. \] 
It is equal to $\max \, \{\alpha_q,\alpha_{q+1}, \dots, \alpha_{e_d-1}\}$. 

\subsection{A summary of results} \label{section.results} 
Define 
\[ \SS(d) = \max \, \{\alpha_1,\alpha_2, \dots, \alpha_{e_d-1} \}, \] 
which is the satiety of $J_2$, and let 
\[ \zeta(d) = \frac{1}{d-2} \sqrt{\frac{(d-1)(d^2-2)}{2}}. \] 
\begin{Theorem} \sl \label{theorem.lubounds} 
For $d \ge 4$, we have inequalities 
\[ \zeta(d) \le \SS(d) \le d+2. \] \end{Theorem} 

Broadly speaking, the lower bound implies that 
$\SS(d)$ grows no slower than $\sqrt{\frac{d}{2}}$. It will be proved in~\S\ref{section.lowerbound}. 
A proof of the upper bound is given in \S\ref{section.spectralseq}. 

The next theorem (which is merely an aggregate of separate propositions) establishes some 
specific lower bounds for $\alpha_1,\alpha_2$ and $\alpha_3$. 
\begin{Theorem} \sl 
Let $(q,b,N)$ denote any of the following triples: 
\[ (1,3,6), \quad (1,4,8), \quad (2,3,12), \quad (3,3,16). \] 
Then $\alpha_q > b$ for all $d \ge N$. 
\label{theorem.triplebounds} \end{Theorem} 
\noindent The proofs are given in~\S\ref{section.triplebounds}. 

\medskip 

The following theorem was inspired by the observation that the saturation sequences 
tend to end in long strings of $3$s. Let 
\begin{equation} 
\begin{array}{llll} 
N_1=4, & N_2=8, & N_3=10, & N_4=14, \\
N_5=18, & N_6=22, & N_7=26, & N_8 = 30. 
\end{array} \end{equation}

\begin{Theorem} \sl 
Let $s$ and $d$ be integers such that $1 \le s \le 8$, and 
$d \ge N_s$. Then at least the last $s$ integers in the saturation sequence of $d$ are all 
equal to $3$. 
\label{theorem.gordan} \end{Theorem} 
\noindent The proof is based upon Gordan's cubic syzygies. It will be given in~\S\ref{section.gordan}. 

\medskip 

In the proofs of the results above, I have had to use machine calculations in order to 
find some complicated compound transvectants, and to evaluate some large determinants. 
They were all done in {\sc Maple}. 

\medskip 

The following two conjectures arise naturally from the previous table. I have 
been unable to make any progress on either of them. 
\begin{Conjecture} \sl \label{conjecture.noni}
The saturation sequence is non-increasing. (This would imply that $\SS(d)=\alpha_1$.) 
\end{Conjecture} 
\begin{Conjecture} \sl 
For all $d \ge 6$, there is always a strict inequality $\alpha_1 > \alpha_2$. 
\end{Conjecture} 

\section{Bounds on $\SS(d)$}  
\subsection{} \label{section.lowerbound} In this section we will prove the lower bound on 
$\SS(d)$. Assume that $(J_2)_m = (I_C)_m$ for some $m >2$. Then the natural morphism 
\[ W_2 \otimes R_{m-2} \lra (I_C)_m \] 
must be surjective, hence by counting dimensions we must have 
\begin{equation} 
(2d-3) \, \binom{m+d-2}{d} \ge \binom{m+d}{d} - (m \, d+1). 
\label{ineq} \end{equation} 
One should like to force a lower bound on $m$ from this inequality. This is carried out in the 
following proposition, which I owe to my colleague A.~Abdesselam. 
Although the proof is elementary in essence, some tricky manipulations are involved. 
\begin{Proposition} \label{prop.lowerbound} \sl 
If $m < \zeta(d)$, then the inequality in~(\ref{ineq}) is false. 
\end{Proposition} 
\demo Transfer the right-hand side of (\ref{ineq}) to the 
left-hand side, and multiply by $d!$. Thus~(\ref{ineq}) is equivalent to 
\[ (2d-3) \left(\prod\limits_{k=m-1}^{m+d-2} k \right) - \left(\prod\limits_{k=m+1}^{m+d} k \right) + 
d! \, (m \, d +1) \ge 0, \] 
or what is the same, 
\begin{equation} 
\underbrace{(2d-3)(m-1)m - (m+d-1) (m+d)}_{Q(d,m)} \times 
\left( \prod\limits_{k=m+1}^{m+d-2} k \right) + d! \, (m \, d+1) \ge 0. 
\label{ineq2} \end{equation} 
We have a factorisation 
\[ Q(d,m) = 2(d-2) ( m - \xi_1(d)) ( m - \xi_2(d)), \] 
where 
\[ \xi_1(d) = \frac{d-1}{d-2} -  \zeta(d), \quad 
\xi_2(d) = \frac{d-1}{d-2} +  \zeta(d). \] 
It is easy to see that $\xi_1(d) <0$ and $\xi_2(d) >0$. 

\medskip 

\noindent{\bf Case $m=3$.} 
After substitution, the left-hand side of~(\ref{ineq2}) becomes 
\begin{equation} 
\begin{aligned} 
- \, & (d^2-7 \, d+24) \, \frac{d! \, (d+1)}{6} + d! \, (3 \, d+1) \\ 
=  - & \frac{d!}{6} \, (d-2) \, (d-2-\sqrt{13}) \, (d-2+\sqrt{13}). 
\end{aligned} \label{expr1} \end{equation} 
Now assume $3 < \zeta(d)$. Then 
\[ 3 (d-2) < \sqrt{\frac{(d-1)(d^2-2)}{2}} < 
\sqrt{\frac{(d-1)(d^2-1)}{2}} = (d-1) \sqrt{\frac{d+1}{2}}, \] 
and since $\frac{d-1}{d-2} \le \frac{3}{2}$ for $d \ge 4$, we have 
\[ 3 < \frac{3}{2} \sqrt{\frac{d+1}{2}}. \] 
This implies that $d >7$, hence (\ref{expr1}) is negative. 

\medskip 

\noindent{\bf Case $m\ge4$.} Assume $m < \zeta(d)$; then 
$\xi_1(d) < 0 < m < \xi_2(d)$, which implies that $Q(d,m) <0$. We want to show that 
left-hand side of (\ref{ineq2}) is negative. Replace 
$m \, d+1$ by the larger quantity $(m+1) \, d$ and divide by $m+1$ to get 
\begin{equation} 
Q(d,m) \times \underbrace{\left(\prod\limits_{k=m+2}^{m+d-2} k\right)}_{T_m} + \, d! \times d. 
\label{expr2} \end{equation} 
It would be sufficient to show that~(\ref{expr2}) is negative. Observe that 
\[ \frac{T_{m+1}}{T_m} = \frac{m+d-1}{m+2} >1, \] 
i.e., $T_m$ increases with $m$. Hence, (\ref{expr2}) is bounded above by 
the quantity 
\begin{equation} 
Q(d,m) \, T_4 + d! \times d = Q(d,m) \, \frac{(d+2)!}{120} + d! \times d . 
\label{expr3} \end{equation} 
Since $m - \xi_1(d) >4$, and $m - \xi_2(d) < \zeta(d) - \xi_2(d) < -1$, we get 
$Q(d,m) < - 8 (d-2)$. Thus (\ref{expr3}) is strictly smaller than 
\[ - 8 \, (d-2) \, \frac{(d+2)!}{120} + d! \times d = 
- \frac{1}{15}(d-4) (d^2+5d+1) \, d! < 0.  \] 
The proposition is proved. \qed 

\subsection{The Koszul complex} \label{section.spectralseq} 
The upper bound on $\SS(d)$ will be established by a spectral sequence argument. 
(Compare the proof of Theorem 1 in~\cite{Shiffman}.)
We refer to~\cite[Ch.~III.5]{Hartshorne} 
for standard results on the cohomology of line bundles on $\P^d$. 
 
The subspace $W_2 \subseteq R_2$ gives a morphism 
\[ S_{2d-4} \otimes \O_{\P^d}(-2) \stackrel{\partial}{\lra} \O_{\P^d}. \] 
By Proposition~\ref{prop.saturation}, we have, $\im \partial = 
\I_C$ (the ideal sheaf of $C$). 
Consider the Koszul complex of $\partial$, and replace $\O_{\P^d}$ with $\I_C$ . This defines a 
complex $\K^\bullet$ of coherent $\O_{\P^d}$-modules 
\[ 0 \ra \K^{-(2d-3)} \ra \dots \ra \K^p \stackrel{h^p}{\ra} \K^{p+1} 
\ra \dots \ra \K^{-1} \stackrel{h^{-1}}{\ra} \K^0 \ra 0, \] 
where 
\[ \K^p = \begin{cases} 
\wedge^{-p} \, S_{2d-4} \otimes \O_{\P^d}(2 \, p) 
& \text{for $-(2d-3) \le p \le -1$,} \\ 
\I_C & \text{for $p=0$.} 
\end{cases} \] 
We will write $\K^\bullet(m)$ for $\K^\bullet \otimes \O_{\P^d}(m)$. 
Let $\HH^p = \text{ker} \, h^{p}/\im h^{p-1}$ denote the cohomology sheaves of 
$\K^\bullet$. 

\subsection{} There are two second quadrant spectral sequences in the range 
\[ -(2d-3) \le p \le 0, \quad 0 \le q \le d,  \] 
which abut to the hypercohomology\footnote{The hypercohomology groups are denoted by 
upper indices on $\H$. There is scarcely any danger of confusion with the covariants $\H_{2q}$, 
which do not appear in this section.}
of $\K^\bullet(m)$; namely 
\begin{equation} \begin{aligned} 
E_2^{p,q} & = H^q(\P^d,\HH^p \otimes \O_{\P^d}(m)), \quad 
\delta_r: E_r^{p,q} \lra E_r^{p-r+1,q+r} \\ 
E_\infty^{p,q} & \Rightarrow \H^{p+q}(\K^\bullet(m)); 
\label{ss1} \end{aligned} \end{equation} 
and 
\begin{equation} \begin{aligned} 
{\widetilde E}_1^{p,q} & = H^q(\P^d,\K^p(m)), \quad 
{\widetilde\delta}_r: {\widetilde E}_r^{p,q} \lra {\widetilde E}_r^{p+r,q-r+1} \\ 
{\widetilde E}_\infty^{p,q} & \Rightarrow \H^{p+q}(\K^\bullet(m)). 
\label{ss2} \end{aligned} \end{equation} 
Henceforth, let 
\begin{equation}  m=d+2. \end{equation}
First, consider the terms in~(\ref{ss1}).  The support of each $\HH^p$ is 
contained in $C$ (see~\cite[Prop.~1.6.5]{BH}), hence 
$E_2^{p,q} = 0$ for $q \ge 2$. This forces $E_2^{p,q} = E_\infty^{p,q}$. 
Since $h^{-1}$ is a surjection, $\HH^0 = 0$. 

The sheaf 
$\HH^{-1}$ will be calculated in Proposition~\ref{prop.H1} below, 
from which it will follow that $H^1(\P^d,\HH^{-1} \otimes \O_{\P^d}(m))=0$. 
Hence $E_2^{p,q}=0$ for all $p+q=0$, implying that 
\begin{equation} \H^0(\K^\bullet(m)) =0. \label{HH0} \end{equation} 

On the other hand, all the nonzero ${\widetilde E}_1^{p,q}$ terms in~(\ref{ss2}) are 
concentrated in the rows $q=0,d$. Our choice of $m$ ensures that 
$\K^{-(d+1)},\K^{-d}$ are respectively equal to 
\[ \wedge^{d+1} S_{2d-4} \otimes \O_{\P^d}(-d), \quad 
\wedge^d S_{2d-4} \otimes \O_{\P^d}(-d+2), \] 
and hence ${\widetilde E}_1^{p,q}=0$ for $(p,q)=(-d-1,d),(-d,d)$. On account of~(\ref{HH0}), this 
forces ${\widetilde E}_\infty^{0,0} = {\widetilde E}_2^{0,0} = 0$. Hence the morphism 
\[ H^0(\K^{-1}(m))  \lra H^0(\K^0(m)) \] 
must be surjective, i.e., $(J_2)_{d+2} = (I_C)_{d+2}$, and thus $\SS(d) \le d+2$. \qed 

\subsection{} \label{section.uv} 
Consider the sheaf $\HH^{-1} = \ker h^{-1}/\im h^{-2}$ supported on $C \simeq \P^1$. 
Henceforth we denote it by $\HH$ for brevity. Since $\K^\bullet$ is 
an $SL_2$-equivariant complex, and the action of $SL_2$ on $C$ is transitive, $\HH$ 
must be torsion-free and hence locally free. 
\begin{Proposition} \sl Assume $d \ge 3$. Then 
$\HH$ is a rank $d-2$ vector bundle on $\P^1$. Moreover, it splits as a direct 
sum of line bundles $\oplus \, \O_{\P^1}(t)$, where each summand satisfies the 
inequalities $-4d+4 \le t \le -2d-2$. 
\label{prop.H1} \end{Proposition} 
\noindent It follows that the group 
\[ H^1(\P^d,\HH \otimes \O_{\P^d}(m)) \simeq 
\oplus \, H^1(\P^1,\O_{\P^1}(m d +t)), \] 
vanishes for $m \ge 4$, since $md+t \ge 4 >-2$. 
This suffices to conclude the argument in the previous section. 

\medskip 

\demo The proof will follow from a calculation of local transition functions. 
Let $\lambda_i = a_i/a_0$ for $1 \le i \le d$, and  
$f = \F/a_0 = x_1^d + \sum\limits_i \, \binom{d}{i} \, \lambda_i \, x_1^{d-i} \, x_2^i$. 
We will write the Hessian $(f,f)_2$ as 
\begin{equation} 
\sum\limits_{r=2}^d \, \binom{2d-4}{r-2} \, u_r \, x_1^{2d-r-2} \, x_2^{r-2} + 
\sum\limits_{s=d+1}^{2d-2}  \, \binom{2d-4}{s-2} \, v_s \, x_1^{2d-s-2} \, x_2^{s-2}; 
\label{affine.hessian} \end{equation} 
where $u_r, v_s$ are elements in the ring $\gA = \QQ[\lambda_1,\dots,\lambda_d]$. 
(The rationale behind this notation will emerge below.) 
A direct calculation with formula~(\ref{formula.transvectants}) 
shows that we have expressions 
\[ \kappa \, u_2 = \lambda_2 - \lambda_1^2, \quad 
\kappa \, u_3 = \lambda_3 - \lambda_1 \, \lambda_2, \] 
and in general 
\[ \kappa \, u_r = \lambda_r - P_r(\lambda_1,\dots,\lambda_{r-1}), \] 
for some polynomials $P_r$. (Throughout, we have used $\kappa$ as a placeholder 
for various nonzero rational constants which need not be precisely specificed. 
See Example~\ref{example.d4} below.) 
If we define the weight of $\lambda_i$ to be $i$, then $u_r,v_s$ are isobaric of weights 
$r,s$ respectively. 

A simple induction shows that 
$\kappa \, u_r \equiv \lambda_r - \lambda_1^r \; \text{mod} \, (u_2,\dots,u_{r-1})$. 
It follows that $u_2,\dots,u_d$ is a regular sequence, and that 
$\ga=(u_2,\dots,u_d) \subseteq \gA$ is 
the defining ideal of the affine piece of $C$ in $\spec \gA \subseteq \P^d$. 

Since $v_s \in \ga$, we must have identities of the form 
$v_s = \sum\limits_{r=2}^d \, g_{s-r} \, u_r$, where $g_{s-r} \in \gA$ are isobaric of 
weight $s-r$. Fix one such an identity for each $s$, and let 
\[ z_s = V_s - \sum\limits_{r=2}^d \; g_{s-r} \, U_r, \quad \text{for} \quad 
d+1 \le s \le 2d-2.\] 

\subsection{} 
Let $M$ denote the free $\gA$-module of rank $2d-3$ on basis elements 
\[ U_r = (-1)^r \, x_2^{2d-r-2} \, x_1^{r-2}, \quad V_s = (-1)^s \, x_2^{2d-s-2} \, x_1^{s-2}, \] 
for the same range of $r,s$ as in~(\ref{affine.hessian}). The notation is chosen in such a way 
that the complex $\K^{-2} \stackrel{h^{-2}}{\lra} \K^{-1} \stackrel{h^{-1}}{\lra} \I_C$ 
is represented over $\spec \gA$ by the $\gA$-module maps 
\[ \wedge^2 \, M \stackrel{\tgf}{\lra} M \stackrel{\gf}{\lra} \ga, \] 
where 
\[ \tgf(W_i \wedge W_j) = w_j \, W_i - w_i \, W_j, \quad \text{and} \quad \gf(W_i) = 
((f,f)_2,W_i)_{2d-4} = w_i. \] 
(Here $W$ stands for either $U$ or $V$ as dictated by the index $i$, and similarly for $w$. 
E.g., $W_2 = U_2, w_{d+1} = v_{d+1}$ etc.) 

\subsection{} Since the $\gA$-module 
\[ N = \Gamma(\spec \gA,\HH) = {\ker \gf}/{\im \tgf}\]  is 
annihilated by $\ga$, it may be regarded as a module over $\gA/\ga \simeq \QQ[\lambda]$. 
(We have written $\lambda$ for $\lambda_1$.) 
It is clear that $z_s \in \ker \gf$. Let $\xi_s$ denote the class of $z_s$ in $N$. 
\begin{Lemma} \sl 
With notation as above, $N$ is the free $\QQ[\lambda]$-module over the elements 
$\{ \xi_s\}$. 
\end{Lemma} 
\demo  If $z = \sum\limits_r \, \alpha_r \, U_r + 
\sum\limits_s \, \beta_s \, V_s \in \ker \gf$, then $z - \sum\limits_s \, \beta_s \, z_s$ is an 
element in $\ker \gf$ which involves only the $U_r$. Hence it must 
necessarily lie in $\im \tgf$, 
since there are no syzygies between the $u_r$ except those coming 
from the tautological Koszul relations. This shows 
that the $\{\xi_s\}$ generate $N$. 
Now consider the map 
\[ e: \QQ[\lambda]^{d-2} \lra N, \quad 
p=(p_{d+1}(\lambda),\dots,p_{2d-2}(\lambda)) \lra 
\sum \, p_s(\lambda) \, \xi_s. \] 
Assume $e(p)=0$, and let $s$ be the largest index such that 
$p_s(\lambda)\neq 0$. Then the weight $s$ part of the relation gives an identity 
$p_s(0) \, \xi_s + \dots = 0$. We may assume that $p_s(0) \neq 0$, since $N$ is torsion-free. 
However, it is clear from the definition of $\tgf$ that no such element can lie in $\im \tgf$. 
Hence $\ker e =0$. \qed 

\subsection{} Now write $\mu_{-i} = a_{d-i}/a_d$ (considered to be of weight $-i$), 
and let $\gA' = \QQ[\mu_{-1},\dots,\mu_{-d}]$. If $f' = \F/a_d$, then $(f',f')_2 =$ 
\[ 
\sum\limits_{r=2}^d \, \binom{2d-4}{r-2} \, u_{-r} \, x_2^{2d-r-2} \, x_1^{r-2} + 
\sum\limits_{s=d+1}^{2d-2}  \, \binom{2d-4}{s-2} \, v_{-s} \, x_2^{2d-s-2} \, x_1^{s-2}; 
\] 
where $u_{-r},v_{-s} \in \gA'$ are isobaric elements of weights $-r,-s$ respectively. 
The same results are true {\sl mutatis mutandis} over $\spec \gA'$, and we have generators 
$\{\xi_{-s}\}$ of $N'$ with weights $-(2d-2),\dots,-(d+1)$. Define the vectors 
\[ \xi^+ = \left[ \begin{array}{c} \xi_{d+1} \\ \vdots \\ \xi_{2d-2} \end{array} \right], 
\quad 
\xi^- = \left[ \begin{array}{c} \xi_{-(2d-2)} \\ \vdots \\ \xi_{-(d+1)} \end{array} \right]. 
\] 
Then $\lambda^{-(3d-1)} \, \xi^+$ and $\xi^-$ are two bases of 
$\Gamma(\spec \gA \cap \spec \gA',\HH)$ as a $\QQ[\lambda,\lambda^{-1}]$-module, 
and hence there is a matrix $Q \in GL(d-2,\QQ[\lambda,\lambda^{-1}])$ such that 
$Q \, \xi^- = \lambda^{-(3d-1)} \, \xi^+$. By taking the weights into account, one sees that 
the $(i,j)$-th entry of $Q$ is of the form $ c \, \lambda^{i-j}$ for some $c \in \QQ$. 

Now apply~\cite[Proposition 3.1]{Haze} to $Q$. It produces a factorisation 
$Q = E^{-1} \, D \,  F$, where 
\[ E \in GL(d-2,\QQ[\lambda]), \quad 
F \in GL(d-2,\QQ[\lambda^{-1}]), \] 
and $D$ is a diagonal matrix of the form 
$\left[ \begin{array}{ccc} \lambda^{k_1} & \dots & 0 \\  & \ddots & \\ 0 & \dots & 
\lambda^{k_{d-2}} \end{array} \right]$. 
Since the entries of 
$\lambda^{d-3} \, Q$ and $\lambda^{-(d-3)} \, Q$ are respectively in $\QQ[\lambda]$ and 
$\QQ[\lambda^{-1}]$, we have $ - (d-3) \le k_i \le d-3$. Hence we have an identity 
\[ F \, \xi^- = 
\left[ \begin{array}{ccc} \lambda^{t_1} & \cdots & 0 \\  & \ddots & \\ 0 & \cdots & 
\lambda^{t_{d-2}} \end{array} \right] \, E \, \xi^+, \] 
where each $t_i$ is sandwiched between $-(3d-1) \pm (d-3)$. This completes the proof of 
Proposition~\ref{prop.H1}. \qed 

\begin{Example} \rm \label{example.d4}
Assume $d=4$, then 
\[ \frac{1}{2} \, u_2 = \lambda_2 - \lambda_1^2, 
\quad u_3 = \lambda_3 - \lambda_1 \, \lambda_2, \quad 
3 \, u_4 = \lambda_4 + 2 \, \lambda_1 \, \lambda_3 - 3 \, \lambda_2^2; \] 
and 
\[ \begin{aligned} 
v_5  & =  3 \, \lambda_1 \, u_4 - 3 \, \lambda_2 \, u_3 + \lambda_3\, u_2, \\
v_6 & = 6 \, \lambda_2 \, u_4 - (2 \, \lambda_3 +6 \, \lambda_1 \, \lambda_2) \, u_3 
+ 3 \, \lambda_2^2 \, u_2 . 
\end{aligned} \] 
Hence 
\[ \begin{aligned} 
\xi_5 & = V_5 -\lambda^3 \, U_2 + 3 \, \lambda^2 \, U_3 - 3 \, \lambda \, U_4, \\ 
\xi_6 & = V_6 - 3 \, \lambda^4 \, U_2 + 8 \, \lambda^3 \, U_3 - 6 \, \lambda^2 \, U_4. 
\end{aligned} \] 
We have an identity 
\[ \underbrace{\left[ \begin{array}{cc} 
-1 & 3 \, \lambda^{-1} \\ 
- 3 \, \lambda & 8 \end{array} \right]}_Q \, 
\left[ \begin{array}{rr} \xi_{-6}' \\ \xi_{-5}' \end{array} \right] = 
\lambda^{-11} \, 
\left[ \begin{array}{rr} \xi_5 \\ \xi_6 \end{array} \right]. \] 
Now, $Q = E^{-1} \, D \, F$ for 
\[ E = \left[ \begin{array}{rr} 3 \, \lambda & -1 \\ -1 & 0 \end{array} \right], \quad 
F = \left[ \begin{array}{cc} 0 & 1 \\ 1 & - 3 \, \lambda^{-1} \end{array} \right], \quad 
D = \left[ \begin{array}{cc} 1 & 0 \\ 0 & 1 \end{array} \right], \] 
and hence 
\[ \left[ \begin{array}{cc}  \xi_{-5}' \\ -3 \, \lambda^{-1} \, \xi_{-5}' +\xi_{-6}' \end{array} 
\right]  = \left[ \begin{array}{cc} \lambda^{-11} & 0 \\ 0 & \lambda^{-11} \end{array} \right] \, 
\left[ \begin{array}{cc} 3 \, \lambda \, \xi_5 - \xi_6 \\ - \xi_5 \end{array} \right], \] 
which gives an isomorphism of $\HH$ with $\O_{\P^1}(-11) \oplus \O_{\P^1}(-11)$. 
\end{Example}

\section{Syzygies in the ring of covariants} 
\subsection{} 
Fix an integer $q$ in the range $1 \le q \le e_d-1$. The following technical result relates 
the magnitude of $\alpha_q$ to the existence of syzygies in the ring $\A$. 
\begin{Lemma} \sl \label{lemma.syzygies} 
For an integer $m \ge 3$, the following conditions are equivalent: 
\begin{enumerate} 
\item[(i)]
$m \ge \alpha_q$. 
\item[(ii)]
Given any covariant $\Phi$ of degree-order $(m-2,n)$, and any integer $r$ such that 
$0 \le r \le \min(2d-4q-4,n)$, there exists an identity of the form 
\begin{equation} 
(\H_{2q+2},\Phi)_r = 
\sum\limits_{i=1}^q \, (\H_{2i}, \Psi_i)_{2(q-i+1)+r+\frac{1}{2}(n_i-n)},  
\label{H2q.syzygy} \end{equation}
for some covariants $\Psi_i$ of degree-orders $(m-2,n_i)$. 
\end{enumerate} 
\end{Lemma}
Broadly speaking, condition (ii) means that any expression of the form 
$(\H_{2q+2},\Box)_\star$ can 
be rewritten as a sum of terms of the form $\{(\H_{2i},\Box)_\star\}_{1 \le i\le q}$ 
using algebraic relations in the ring $\A$. The index of transvection of the term 
$(\H_{2i},\Psi_i)$ is determined by the requirement that each summand should have order 
$2d-4q-4+n-2r$ in $\ux$. 
\begin{Example} \rm 
Assume $d=4$, and let $(q,m)=(1,3)$. The only choice for $\Phi$ (up to a constant) is $\F$, 
and since $\H_4$ is an invariant, $r=0$. We have an identity $\H_4 \, \F = 6 \, (\H_2,\F)_2$ 
(see~\cite[\S 93]{GrYo}), hence condition (ii) is satisfied. This shows that $\alpha_1 = 3$. 
\end{Example} 

\begin{Example} \rm Assume $d=7$. The space $\A_{3,9}$ is two dimensional, and it is easy to 
show (say by specialising $\F$) that $\{(\H_4,\F)_2, (\H_2,\F)_4 \}$ is a basis. 
Hence there is no identity of the type (\ref{H2q.syzygy}) for $(q,m,r)=(1,3,2)$ and $\Phi = \F$, 
which shows that $\alpha_1 > 3$. 

On the other hand, if one takes $(q,m)=(1,4)$, then such identities always exist. For instance, 
if $\Phi = \H_6$ and $r=2$, then 
\[ (\H_4, \H_6)_2 = \frac{42}{13} \, (\H_2, \F^2)_{10} + \frac{15876}{845} \, 
(\H_2, \H_2)_8 + \frac{10332}{715} \, (\H_2, \H_4)_6. \] 
This can be verified by the use of symbolic calculus as in~\cite[Ch.V]{GrYo}.
\end{Example} 

\medskip 

\noindent {\sc Proof of Lemma~\ref{lemma.syzygies}.}
Let $\U$ denote the image of the morphism 
\[ W_{2q+2} \otimes R_{m-2} \lra R_m. \] 
By definition, it is spanned by all the coefficients of all the transvectants of the form $(\H_{2q+2},\Phi)_r$. 
Similarly $(J_{2q})_m$ is spanned by the union of images of the maps 
\[ W_{2i} \otimes R_{m-2} \lra R_m, \quad (1 \le i \le q).  \] 
The inequality $m \ge \alpha_q$ holds iff $\U$ is contained in $(J_{2q})_m$, which happens iff 
an arbitrary $(\H_{2q+2},\Phi)_r$ can be rewritten as in~(\ref{H2q.syzygy}). 
This proves the lemma. \qed 

\subsection{} For what it is worth, the lemma gives some 
thematic support to Conjecture~\ref{conjecture.noni}. Indeed, as $m$ is held constant 
and $q$ decreases, the 
range of allowable values of $r$ increases and hence, {\sl prima facie}, condition (ii) 
becomes more stringent. This makes it plausible that $\alpha_q$ should increase (or at 
least remain stationary) with decreasing $q$. 

\subsection{} \label{section.triplebounds} 
The next four propositions are the ingredients in Theorem~\ref{theorem.triplebounds}. 
In each case we establish a lower bound on some $\alpha_q$ by showing that a certain type of 
syzygy cannot exist in $\A$ for sufficiently large $d$. 

\begin{Proposition} \sl 
If $d \ge 12$, then $\alpha_2 >3$. 
\end{Proposition} 
\demo Let $(q,m)=(2,3), \Phi = \F$, and $r = 6$ in the notation of 
Lemma~\ref{lemma.syzygies}. To show that condition (ii) fails, it 
is enough to show that the set 
\[ \Gamma_1 = (\H_6,\F)_6, \quad \Gamma_2 = 
(\H_4,\F)_8, \quad \Gamma_3 = (\H_2,\F)_{10},  \] 
is linearly independent. Specialise to the form 
\[ F = x_1^d + x_1^{d-2} \, x_2^2 + x_1 \, x_2^{d-1} + x_2^d, \] 
and calculate the $\Gamma_i$. Construct a $3 \times 3$ matrix $M$ whose $i$-th row 
sequentially consists of the coefficients of 
\[ x_1^{2d-12}\, x_2^{d-12}, \quad x_1^{2d-13}\, x_2^{d-11}, \quad x_1^{2d-15}\, x_2^{d-9} \] 
in $\Gamma_i$. For instance, the $(2,1)$-entry is 
\[ \frac
{(d-8) \, (d-9) \, (d-10) \, (d-11)}
{8 \, (2 \, d-9) \, (2 \, d -11) \, (2 \, d - 13) \, (2 \, d - 15)}. \] 
Now $\det(M)$ is a rational function in $d$, and one easily checks (in {\sc Maple}) that it is 
nonzero for $d \ge 12$. \qed 

\medskip 

One needs to expend a certain quantity of trial and error to discover that $r=6$ would make the 
proof work. The analogous argument fails for the set 
\[ (\H_6,\F)_r, \quad (\H_4,\F)_{r+2}, \quad (\H_2,\F)_{r+4}, \] 
if $r=0,1,2,3,4,5$. Similar remarks apply to the results below. 

\medskip 

\begin{Proposition} \sl If $d \ge 6$, then $\alpha_1 >3$. 
\end{Proposition} 
\demo It is enough to show that $(\H_4,\F)_2$ is not a constant multiple of $(\H_2,\F)_4$
for $d \ge 6$. This is done by specialising to the same $F$ as above. \qed 

\begin{Proposition} \sl 
If $d \ge 16$, then $\alpha_3 >3$. 
\end{Proposition} 
It is enough to show that $(\H_8,\F)_{10}$ cannot be written as a linear combination of 
\begin{equation} 
(\H_6,\F)_{12}, \quad (\H_4,\F)_{14}, \quad (\H_2,\F)_{16}, 
\label{list.alpha3} \end{equation} 
which can be checked by specialising to 
$F = x_1^d + x_1^{d-3} \, x_2^3  - x_1 \, x_2^{d-1} + 2 \, x_2^d$. 
The details are similar to above. However, this argument works only for $d \ge 18$. 
If $d=16,17$, then unfortunately $(\H_8,\F)_{10}$ is linearly dependent on the 
three covariants in (\ref{list.alpha3}), hence one has to look for specific features of those 
cases. 

Assume $d=16$ or $17$, and let $\Phi = \F$ and $(q,r)=(3,16)$. 
One can check by specialisation that the 
covariant $(\H_8,\F)_{16}$ does not vanish identically for $d=16,17$. 
It is clear that no relation of the type (\ref{H2q.syzygy}) can exist, since the 
index of transvection in each summand on the right must be at least $18$, which is impossible. 
This completes the proof. \qed

\begin{Proposition} \sl For $d \ge 8$, we have $\alpha_1 > 4$. 
\end{Proposition} 
\demo 
It is enough to show that there is no constant $\eta_d \in \QQ$ such that 
\[ \underbrace{(\H_4,\H_4)_{2d-8}}_J = \eta_d \, 
\underbrace{(\H_2,\H_2)_{2d-4}}_K. \] 
Let 
\[ F_1 = x_1^d + x_2^d, \quad F_2 = x_1^d + x_1^{d-2} \, x_2^2 + x_1\, x_2^{d-1}, \] 
and consider the determinant 
$ \left| \begin{array}{cc}  
J_1 & J_2 \\ K_1 & K_2 \end{array} \right|$, where $J_i,K_i$ denote the specialisations 
of those invariants to $F_i$. It is enough to show that this determinant does not vanish 
for any $d \ge 8$. 
An explicit calculation shows that up to a nonzero factor, it equals 
\[ f(d) = \underbrace{(d^3 - 8 \, d^2 + 19 \, d - 14)}_{T_1} + \, (-1)^d \, 
\underbrace{\binom{2d-6}{d-3}}_{T_2}. \] 
There is nothing to show for even $d$, so assume it to be odd. 
Now $d^4 >T_1$ (because $d^4-T_1$ has no real 
roots) and 
$T_2 > 2^{d-3}$. For $d \ge 21$, we have $2^{d-3} > d^4$, and hence $f(d) \neq 0$. Thus 
it only remains to verify the claim for $d = 9,11,\dots,19$, which is routine. \qed 

\medskip 

In general, let $G^{(q)} = (\H_{2q},\H_{2q})_{2d-4q}$, which is a degree $4$ invariant of 
$d$-ics, moreover the $\{G^{(q)}\}$ span the space $\A_{4,0}$. 
One can deduce a formula for the number $h(d) = \dim \A_{4,0}$ as follows. 
By Hermite reciprocity (see~\cite[\S 157]{Salmon}), 
it is the same as the number of linearly independent invariants of degree $d$ for binary quartics. 
If $\F$ denotes the generic quartic, then each such invariant is necessarily of the form 
$[(\F,\F)_4]^a \, [(\F,(\F,\F)_2)_4]^b$. 
Hence $h(d)$ is the cardinality of the set  
\[ \{(a,b) \in {\mathbf N}^2: 2 \, a + 3 \, b = d \}. \] 
This gives the following formula: write $d = 6 \, e + k$ where $0 \le k \le 5$. Then 
$h(d) = e+\delta_k$, where $\delta_1 = 0$ and $\delta_k = 1$ for $k \neq 1$. 
For instance, $h(75) = 13$. 
\begin{Proposition} \sl 
In the saturation sequence of $d$, at least $h(d)$ of the integers are strictly greater than $4$. 
\end{Proposition} 
\demo Assume that $G^{(q_i)}, (i = 1,2, \dots, h)$ are linearly independent. Then it is 
immediate that each $\alpha_{q_i} > 4$. \qed 

\medskip 

The results in this section, a little scattered and unsystematic as they are, should be illustrative 
of the principle that in so far as the syzygies in $\A$ are intricate and unruly 
(e.g., see~\cite[Ch.~VII]{GrYo} or~\cite{Bedratyuk}), 
it seems unlikely that one can deduce precise formulae for the $\alpha_q$. 
\section{Gordan's syzygies} \label{section.gordan} 
We begin with an explanation of Gordan's cubic syzygies (see~\cite[\S 54]{GrYo}). 
They will be used to prove Theorem~\ref{theorem.gordan}. 

\subsection{} 
Let $f,\phi,\psi$ denote binary forms of orders $m,n,p$ respectively; and let $a_1,a_2,a_3$ be 
nonnegative integers such that 
\[ a_2 + a_3 \le m, \quad a_1 + a_3 \le n, \quad a_1 + a_2 \le p. \] 
Assume furthermore, that at least one of the following conditions is true: 
\[ a_1 = 0, \quad \text{or} \quad a_2 + a_3 = m. \] 
Then Gordan's syzygy (or series) is the identity 
\[ \begin{aligned} 
{} & \sum\limits_{i =0}^\infty \;  
\frac{\binom{n-a_1-a_3}{i} \binom{a_2}{i}}{\binom{m+n-2a_3-i+1}{i}}
\, ((f,\phi)_{a_3+i},\psi)_{a_1+a_2-i} \\ 
= \, (-1)^{a_1} \, & 
\sum\limits_{i =0}^\infty \; 
\frac{\binom{p-a_1-a_2}{i}\binom{a_3}{i}}{\binom{m+p-2a_2-i+1}{i}} 
\, ((f,\psi)_{a_2+i},\phi)_{a_1+a_3-i}. \end{aligned} \] 
It is usually denoted by $\left( \begin{array}{ccc} f & \phi & \psi \\ m & n & p \\ a_1 & a_2 & a_3 
\end{array} \right)$. 
By convention, $\binom{a}{b} =0$ if $a<b$, hence either side is a finite sum. The 
total index of transvection in each term is $a_1+a_2+a_3$, which is also called the 
{\sl weight} of the syzygy. In the following 
two sections we will specialise to the case $f = \phi = \psi = \F$, and rewrite the syzygies in 
a more convenient form. 

Let $\{a,b\}$ denote the cubic covariant $((\F,\F)_a,\F)_b$ of order $3d-2(a+b)$. It vanishes 
identically unless 
\begin{equation} 
0 \le a, \, b \le d, \quad \text{$a$ is even \; and} \quad 2a + b \le 2d. 
\label{adm.pair} \end{equation} 
An {\sl admissible pair} $(a,b)$ is one which satisfies the conditions in~(\ref{adm.pair}). 
(However, these conditions do not guarantee that $\{a,b\}$ is nonzero; e.g., 
if $d=5$, then $\{2,5\}$ vanishes identically -- see~\cite[\S 71]{GrYo}.) 

\subsection{Syzygies of weight at most $d$.} Choose integers $w,k$ in the range 
\[ 0 \le w \le d, \quad 0 \le k < \frac{w}{2}, \] 
and let $a_1=0, \, a_2=k, \, a_3=w-k$. 
Then we have a syzygy 
\begin{equation} \Gordan_\bullet(k,w): \; \sum\limits_{m = k}^w \; 
\theta_{d,k,w}^{(m)} \, \{m,w-m \} = 0,  \label{Gordan1} \end{equation}
where 
\[ \theta_{d,k,w}^{(m)} = 
\frac{\binom{d-k}{m-k} \binom{w-k}{m-k}}{
\binom{2d-k-m+1}{m-k}} - 
\underbrace{\frac{\binom{d-w+k}{m-w+k} \binom{k}{m-w+k}}{
\binom{2d-w+k-m+1}{m-w+k}}}_{(\star)}.  \] 
The term $(\star)$ is understood to be zero if $m <w-k$. 
For instance, if $d=7$, then $\Gordan_\bullet(1,6)$ is the syzygy 
\[ \frac{5}{2} \, \{2,4\} + \frac{5}{3} \, \{4,2\} - \frac{11}{28} \, \{6,0\}=0. \] 
\subsection{Syzygies of weight at least $d$.} 
Alternately, choose integers $w,k$ in the range 
\[ d \le w \le \frac{3 \, d}{2}, \quad w-d \le k \le \frac{d}{2}, \] 
and let $a_1 = w-d, \; a_2 = d-k, \; a_3 = k$. 
Then we have a syzygy 
\begin{equation} 
\Gordan^\bullet(k,w) : \sum\limits_{m =k}^{2d-w} \; 
\vartheta_{d,k,w}^{(m)} \, \{m,w-m\} = 0, \label{Gordan2} \end{equation} 
where 
\[ \vartheta_{d,k,w}^{(m)} = 
\frac{\binom{2d-w-k}{m-k} \binom{d-k}{m-k}}{\binom{2d-k-m+1}{m-k}} + 
(-1)^{w+d+1} \, 
\underbrace{\frac{\binom{d-w+k}{m-d+k} \binom{k}{m-d+k}}
{\binom{d-m+k+1}{m-d+k}}}_{(\star\star)}. \] 
The term $(\star \star)$ is understood to be zero if $m <d-k$. 
For instance, if $d=11$, then $\Gordan^\bullet(4,13)$ is the syzygy 
\[ \{4,9\} + \frac{35}{13} \, \{6,7\} - \frac{31}{66} \, \{8,5\} = 0. \] 
The syzygies $\Gordan^\bullet(k,d)$ and $\Gordan_\bullet(k,d)$ are identical. 

\subsection{} Let us prove Theorem~\ref{theorem.gordan} for $s=1$, which claims that 
$\alpha_{e_d-1}$ is always equal to $3$. First, assume $d$ is even, 
then $\H_d$ is an invariant. It is sufficient to 
show the existence of a syzygy~(\ref{H2q.syzygy}) for $\Phi = \F$ and $r=0$. 
This follows from the fact that the coefficient of $\{d,0\}$ in $\Gordan_\bullet(1,d)$ is 
\[ \theta^{(d)}_{d,1,d} = \frac{1}{d} - \frac{1}{2} \neq 0. \] 

If $d$ is odd, consider the coefficients of $\{d-1,0\}, \{d-1,1\}, \{d-1,2\}$ in 
the syzygies $\Gordan_\bullet(1,d-1), \Gordan_\bullet(1,d)$ and $\Gordan^\bullet(1,d+1)$ 
respectively. They are 
\[ \frac{6}{d(d+1)}-\frac{1}{2}, \quad \frac{6 \, (d-1)}{d \, (d+1)} - 1, \quad 
\frac{6}{d(d+1)}+1, \] 
none of which can be zero. This completes the argument. \qed 

\subsection{} 
The following example should illustrate the idea behind the 
proof of Theorem~\ref{theorem.gordan}. 
Suppose we want to show that $\alpha_2=3$ for $d=9$. This requires 
showing (amongst other things) that $\{6,2\}$ can be written as a linear combination of 
$\{2,6\}$ and $\{4,4\}$. However, any of the Gordan syzygies involving 
$\{6,2\}$ will also involve the unwanted term $\{8,0\}$. 
One can use two syzygies simultaneously in order to 
eliminate the latter. For instance, $\Gordan_\bullet(1,8)$ and 
$\Gordan_\bullet(2,8)$ can be written as 
\[ \begin{aligned} 
\frac{49}{33} \{6,2\} -\frac{13}{30} \, \{8,0\} & = 
- \frac{7}{2} \{2,6\} - \frac{70}{13} \, \{4,4\}, \\ 
\frac{13}{22} \{6,2\} -\frac{13}{60} \, \{8,0\} & = 
- \{2,6\} - \frac{105}{26} \, \{4,4\}. 
\end{aligned} \] 
Since the determinant $\left| \begin{array}{rr} 
{49}/{33} & -{13}/{30} \\ 
{13}/{22} & -{13}/{60} \end{array} \right|$ is nonzero, 
$\{6,2\}$ is expressible as a linear combination of $\{2,6\}$ and $\{4,4\}$. The argument in 
the general case is conceptually the same, but the technical details are somewhat tedious. 

\subsection{} Given an admissible pair $(a,b)$, define its {\sl position} $p(a,b)$ to be the number of 
admissible pairs $(a',b')$ of the same weight 
such that $a \le a'$. In any Gordan syzygy involving $\{a,b\}$, it is the $p(a,b)$-th term from the right. 
For instance, if $d=13$, then the sequence $(6,9),(8,7),(10,5)$ shows that $p(6,9)=3$. 

Fix a positive integer $s$. Our object is to find an integer $N_s$ such that $\alpha_{e_d-s} =3$ for 
$d \ge N_s$. We will assume that $d \ge 4s-2$; this will prove useful in manipulating 
the syzygies. (We are making no attempt to find the optimal value of $N_s$.) 
First, assume $d$ to be even, say $d = 2n$. Then 
$\H_{2(n-s+1)}$ has order $4s-4$, and hence the possible candidates for the left-hand side 
of (\ref{H2q.syzygy}) are 
\[ \{2 \, (n-s+1),t\}, \quad \text{for $0 \le t \le \min(d,4s-4) = 4s-4$}. \] 
Let $w = 2 \, (n-s+1)+t$. 

\medskip 

\noindent {\bf Case I. } Assume $0 \le t \le 2 s-2$, then $w \le d$. 
It is easy to see that the position $p = p(2n-2s+2,t)$ equals $[\frac{t}{2}]+1$. 
Construct a $p \times p$ matrix $M_t$ whose $(k,m)$-th 
entry is $\theta^{(2m)}_{d,k,w}$, for 
\[ 1 \le k \le p, \qquad n-s+1 \le m \le n-s+p. \] 

\noindent {\bf Case II. } Assume $2s-1 \le t \le 4s-4$, then 
$d+1 \le w \le \frac{3d}{2}$ and $p = 2s-1-\lceil \frac{t}{2} \rceil$. Construct $M_t$ by 
letting its $(k,m)$-th element to be $\vartheta^{(2m)}_{d,k,w}$, for 
\[ w-d \le k \le w-d+p-1, \qquad n-s+1 \le m \le n-s+p. \] 

\medskip 

Now let $d$ be odd, say $d = 2n+1$. Then $\H_{2(n-s+1)}$ has order $4s-2$, and one can construct 
matrices $M_t'$ as above for $0 \le t \le 4s-2$. It is clear that 
\[ \Delta_t(d) = \det M_t, \quad \Delta_t'(d) = \det M_t', \] 
are rational functions of $d$. I have calculated them explicitly for $s \le 8$, and in each case 
determined the threshold $N_s$ such that they are all nonzero for $d \ge N_s$. 
The computations were programmed in {\sc Maple}. 
For instance\footnote{It seems to be a general feature that the numerators and 
denominators of $\Delta_t,\Delta_t'$ almost entirely consist of linear factors. Why this should be
so is not obvious to me.}, if $s=3$, then 
\[ \Delta_6(d) = 
\frac{3780 \, (d-4) \, (d-5) \, (d-6) \, (d+7) \, (d^2+3 \, d + 10)}
{(d-1)^2 \, (d-2) \, (d+2) \, (d+1)^2 \, d^2 \, (d+3)}, \] 
which is nonzero for $d >6$. 

As in the example 
above, this shows the existence of a syzygy for each $\{2(n-s+1),t\}$ 
as required by~(\ref{H2q.syzygy}). \qed 

\medskip 

The argument would break down if any of the determinants were to vanish identically; 
but fortunately this does not happen, at least for $s \le 8$. The theorem could 
be mechanically extended to a few more values of $s$, but this is unlikely to be of much 
interest in itself. This line of argument suggests the following conjecture. 

\begin{Conjecture} \sl For any positive integer $s$, there exists an integer $N_s$ such that 
(at least) the last $s$ integers in the saturation sequence are equal to $3$ for all $d \ge N_s$. 
\label{conjecture.tail3} \end{Conjecture} 

This would follow immediately if it could be shown that $\Delta_t,\Delta'_t$ never vanish 
identically. Furthermore, the data suggest that $N_s = 4 \, s -2$ is in fact the best possible value for 
$s \ge 3$. 

\bigskip 

{\small 
\noindent {\sc Acknowledgements:} I thank Isabella Bermejo and Vijay Kodiyalam for some 
instructive correspondence, and Abdelmalek Abdesselam for 
Proposition~\ref{prop.lowerbound}. I am also grateful to Daniel Grayson and 
Mike Stillman (the authors of Macaulay-2). This work was financially supported in part by 
NSERC, Canada.}

\medskip 

\centerline{--} 

\vspace{1cm}

\parbox{7cm}{ \small 
Jaydeep Chipalkatti \\
Department of Mathematics \\ 
University of Manitoba \\ 
Winnipeg, MB R3T 2N2 \\ 
Canada. \\ 
{\tt chipalka@cc.umanitoba.ca}}

\end{document}